# Recognizing badly presented $Z$-modules.

George Havas [*]
Key Centre for Software Technology
Department of Computer Science
University of Queensland
Queensland 4072 Australia

Derek F. Holt
Mathematics Institute
University of Warwick
Coventry CV4 7AL
UK

Sarah Rees
Department of Mathematics and Statistics
University of Newcastle
Newcastle upon Tyne NE1 7RU
UK

**Abstract**

Finitely generated $Z$-modules have canonical decompositions. When such modules are given in a finitely presented form there is a classical algorithm for computing a canonical decomposition. This is the algorithm for computing the Smith normal form of an integer matrix. We discuss algorithms for Smith normal form computation, and present practical algorithms which give excellent performance for modules arising from badly presented abelian groups.

We investigate such issues as congruential techniques, sparsity considerations, pivoting strategies for Gauss-Jordan elimination, lattice basis reduction and computational complexity. Our results, which are primarily empirical, show dramatically improved performance on previous methods.

## 1 Introduction

Hartley and Hawkes [13, Chapter 8], describe canonical decompositions for finitely generated modules over principal ideal domains. In Chapter 10, they refine this to finitely generated abelian groups and their associated $Z$-modules. This is based on theorems about matrices over principal ideal domains and their canonical forms in Chapter 7, where they give an algorithm for the computation of a canonical form. For integer matrices ($Z$-modules) these results date back to Smith [33].

Our motivation is the solution of group-theoretic problems, so we cast our description in terms of abelian groups. Since $Z$-modules are no more nor less than abelian groups, the principles are relevant for $Z$-modules in general. Sims [32] includes a significant chapter on abelian groups in his book on computational group theory.

[*]Partially supported by Australian Research Council Grant A49030651.

A finitely presented abelian group $G$ may be given by a set of $n$ generators $x_1, \ldots, x_n$ and $m$ relations $\sum_{j=1}^{n} a_{i,j} x_j = 0$. Such presentations arise in various types of computation in natural ways. Examples include subgroup presentation by Reidemeister-Schreier processes, cohomology calculations, and as part of soluble quotient computation algorithms (see Cannon and Havas [4] for an overview and references).

Our aim is to identify such a finitely presented $G$ using effective algorithms. In principle it is easy, but it can be difficult if $n$ and $m$ are large. The fundamental theorem for finitely generated abelian groups tells us all about such groups.

**Theorem.** Let $G$ be a finitely-generated abelian group. Then $G$ has a direct decomposition

$$G \cong G_1 \oplus \ldots \oplus G_r \oplus G_{r+1} \oplus \ldots \oplus G_{r+f}$$

where:
1) $G_i$ is a nontrivial finite cyclic group of order $l_i$ for $i = 1, \ldots, r$;
2) $G_i$ is an infinite cyclic group for $i = r+1, \ldots, r+f$;
3) $l_1 \mid l_2 \mid \ldots \mid l_r$.

The integers $f$ and $l_1, \ldots, l_r$ occurring in such a decomposition are uniquely determined by $G$.

To recognize $G$ we need to determine these integers. Sometimes we want even more. We may also want to know an isomorphism between $G$ as originally given and its canonical form. There is a "standard" algorithm for finitely presented abelian groups, based on matrix diagonalization using Gauss-Jordan elimination. Unfortunately the standard algorithm suffers from serious practical difficulties.

We associate with $G$ its relation matrix, the $m \times n$ integer matrix $A = (a_{i,j})$. $A$ and $B$ are equivalent if there exist unimodular $P$ and $Q$ such that $PAQ = B$. Abelian groups with equivalent relation matrices are isomorphic. An arbitrary integer matrix is equivalent to a unique matrix of the form

$$\begin{pmatrix} b_{1,1} & 0 & 0 & \ldots & 0 \\ 0 & b_{2,2} & 0 & \ldots & 0 \\ 0 & 0 & \ddots & & \\ 0 & \ldots & & & \end{pmatrix}$$

where the $b_{i,i}$ are nonnegative integers and $b_{i-1,i-1} \mid b_{i,i}$. This is the Smith normal form.

Equivalence of matrices can be characterized in terms of elementary row and column operations. For integer matrices the elementary operations are:
- multiply a row (column) by $-1$;
- interchange two rows (columns);
- add an integer multiple of one row (column) to another.

An algorithm for computing the SNF follows (essentially the algorithm given in Hartley and Hawkes).

The first stage of the reduction is to compute an equivalent form

$$\left( \begin{array}{c|ccc} d_1 & 0 & \ldots & 0 \\ \hline 0 & & & \\ \vdots & & B_1 & \\ 0 & & & \end{array} \right)$$

where $d_1$ divides every entry in the submatrix $B_1$.

If $A$ is the zero matrix we are finished. If not, choose a nonzero entry (the *pivot*) and move it to $a_{1,1}$ by suitable row and column interchanges.

While there is an entry $a_{1,j}$ in the first row not divisible by $a_{1,1}$: compute $a_{1,j} = a_{1,1}q + r$ using the Euclidean algorithm; subtract $q$ times the first column from the $j$th column; interchange the first and $j$th columns. Do the same with rows and columns interchanged.

After this, $a_{1,1}$ divides every entry in its row and column. Subtract suitable multiples of the first column (row) from the other columns (rows) to replace each entry in the first row (column), except $a_{1,1}$, by zero. Then we have the correct shape. If the divisibility condition is satisfied, we have finished. If not, there is an entry $a_{i,j}$ such that $a_{1,1}$ does not divide $a_{i,j}$: then add the $i$th row to the first row and return to the while statement. (In practice many algorithms simply compute a diagonal form first, sorting out the divisibility along the diagonal as a final step. In such algorithms this paragraph is delayed to the end.)

Finally reduce $B_i$ recursively. (Efficient implementations of this algorithm use iteration rather than recursion.)

Thus the Smith normal form (in group theoretic terms) leads to methods for decomposing a finitely presented abelian group into a canonical direct product of cyclic subgroups. The $b_{i,i}$ are called the elementary factors of $A$, while the greatest $r$ such that $b_{r,r} \neq 0$ is the rank of $A$. The nontrivial elementary factors of $A$ are the torsion invariants of the associated group $G$ and $n - r$ is the torsionfree rank.

Smith described how to compute the elementary factors in terms of gcd's of subdeterminants of $A$. Thus $b_{i,i}$ is the gcd of all the $i \times i$ subdeterminants of $A$ divided by the gcd of all the $(i-1) \times (i-1)$ subdeterminants of $A$. This description is no good for practical computation with large matrices. It has exponential complexity in terms of $n$, $m$. Furthermore it does not give an isomorphism.

The presentations for $Z$-modules which arise in our area of interest are bad in the sense that they are very distant from the canonical presentation given by the Smith normal form. Thus, the abelian group presentations arising from a Reidemeister-Schreier process (see Havas [14], Havas and Sterling [18] and Neubüser [27]) have large numbers of generators and relations, with a substantial number of trivial modules in the canonical form.

We want efficient algorithms for recognizing the badly presented $Z$-modules which arise, and sometimes also for determining an isomorphism between the initial presentation and a canonical presentation. This can be done with good algorithms for Smith normal form computation for the associated integer matrices.

The problem and some solutions for the group-theoretic context are studied in detail by Havas and Sterling [18] and Sims [32]. Many other researchers have investigated this problem in various contexts. Useful references not included by Havas and Sterling are Hu [20, Appendix A], Gerstein [11], Frumkin [9, 10], Kannan and Bachem [23], Chou and Collins [5], Domich, Kannan and Trotter [7], Domich [6], Iliopoulos [21, 22], Donald and Chang [8] and Hafner and McCurley [12].

In this paper we discuss both integer based methods and modular techniques. We first describe the relevant methods and then present examples of their performance.

## 2 Gauss-Jordan elimination

Superficially methods based on Gauss-Jordan elimination over $Z$ look attractive. Naively implemented the complexity is polynomial, $O(X^3)$ where $X = \max(n, m)$. It does give an isomorphism, which is specified by the unimodular matrix $Q$ and readily computed in the process. But something very bad happens. During Gauss-Jordan elimination over $Z$ the entries in the matrix are easily seen to be bounded **exponentially** in length. Frumkin [9] gives the easy to obtain upper bound of $x^{3^k}$ at the $k$th step (where $x = \max |a_{i,j}|$).

If entries do increase exponentially in length then what appeared to be $O(X^3)$ complexity turns into a polynomial number of arithmetic operations on exponentially large operands. However nobody has shown that any variant of Gauss-Jordan elimination does actually lead to exponential growth.

In fact polynomial bounds on operand size have been found for variants of Gauss-Jordan elimination over $Z$. Polynomial bounds are given by both Kannan and Bachem [23] and improved by Chou and Collins [5] for specific implementations of elimination methods. Our examples show that the Kannan-Bachem and Chou-Collins methods are substantially worse than ours for the kinds of $Z$-modules in which we are interested.

In practice "entry explosion" often occurs when Gauss-Jordan elimination is performed over $Z$. Naive implementations which try to solve this problem using multiple precision arithmetic take too long, while other implementations just blow up.

We address this difficulty by using heuristic pivot selection strategies and some lattice basis reduction methods. Alternatively, various people have observed that modular techniques can be used to avoid entry explosion.

## 3 Modular techniques

Congruential techniques can be used to compute matrix ranks and determinants. In such cases these techniques are fast and avoid entry explosion, with calculations being done in prime fields, $Z_p$ instead of $Z$ (see, for example, Cabay and Lam [2]). For Smith normal form calculation the situation is somewhat more complicated. Underlying theory for modular techniques may be found in Gerstein [11]. Perhaps the first algorithm description, though incomplete, is in Hu [20, Appendix A]. This type of method has been described by Frumkin [9, 10], Havas and Sterling [18], Domich, Kannan and Trotter [7], Domich [6] and Iliopoulos [21, 22], The following description is the basis for the algorithm implemented by Havas and Sterling.

If $G \cong T \oplus F$ where $T$ is torsion and $F$ torsionfree, then $|T| = \prod_1^r b_{i,i}$. By Smith, $|T|$ = the gcd of all the $r \times r$ subdeterminants of $A$. This still does not look promising because the number of subdeterminants is exponential in $n$.

However a useful multiple of $|T|$ may come from the gcd of a small number of subdeterminants. So we get the algorithm outline: determine $r$; calculate $S$, a multiple of $|T|$; perform Gauss-Jordan elimination in $Z_S$, the ring of integers modulo $S$.

Consider $G \cong T \oplus F$. We have $n$ generators and $m$ relators for $G$ and an $m \times n$ relation matrix $A$. We want the SNF for $A$ plus (possibly) an isomorphism $G \to T \oplus F$.

Start by finding the rank $r$ of $T$ and rank $f (= n - r)$ of $F$. This can be done by Gauss-Jordan eliminations over $Z_p$ for a number of primes $p$.

In Havas and Sterling $r$ is "guessed" from one computation modulo a "random" large prime. We have never seen this guess to be wrong in practice, though it is possible to construct examples where it fails. If you are happy to use the guess, or alternatively the guess shows $A$ to have full rank, which must be right, then this step is $O(X^3)$.

To ensure correct rank computation we can use a bound on determinants, e.g., the Hadamard bound, or more crudely, $R^X$, where $R$ is the largest row norm in the matrix. We compute the rank modulo a set of distinct primes whose product exceeds a bound on any subdeterminant. It follows that at least one of the primes does not divide $|T|$. The rank modulo that prime (the maximum of the ranks modulo each of the primes) is the correct rank.

A detailed complexity analysis of the method appears in Hafner and McCurley [12]. They show that it has bit complexity $O(X^4 W(X \log xX) \log xX)$, where $W$ is a function which measures the number of bit operations required to compute with numbers the size of its argument. $W(y)$ is $O(y \log^2 y \log \log y)$. This gives a bit complexity which is $\tilde{O}(X^5)$, ignoring logarithmic terms. An informal analysis follows, based on a RAM model of computation.

Assume the primes we use have size $P$. Then this computation is $O(X^3 . \log(R^X)/\log(P))$ at worst, i.e. $O(X^4 . \log(R)/\log(P))$. (The constants here are good, since it is just Gauss-Jordan elimination over $Z_p$ and $\log(R)/\log(P)$ may be significantly less than 1.)

If $F$ is trivial, or if we do not want an isomorphism, then try to find a suitable multiple $S$ of $|T|$ by computing the gcd of a few $r \times r$ subdeterminants. That is not worse than another $O(X^4 . \log(R)/\log(P))$ calculation. (Note that rank and determinant computations amount to essentially the same calculation.)

Finally do a Smith normal form calculation in $Z_S$, which is an $O(X^3)$ computation as long as $S$ is "small". If $F$ is trivial we readily get an isomorphism $G \to T$. All in all the above behaves like an $O(X^3)$ calculation, in spite of being formally $O(X^4)$. The $\log(R^X)/\log(P)$ term, which contributes one power of $X$, specifies how many times Gauss-Jordan eliminations need to be done, and is more like a constant.

Even if $S$ is **bad** this approach can be successfully extended using a different decomposition.
**Theorem.** Let $G$ be a finitely-generated abelian group. Then $G$ has a direct decomposition

$$G \cong G_1 \oplus \ldots \oplus G_s \oplus G_{s+1} \oplus \ldots \oplus G_{s+f}$$

where:
1) $G_i$ is a nontrivial finite cyclic group of prime-power order $p_i^{\alpha_i}$ for $i = 1, \ldots, s$;
2) $G_i$ is an infinite cyclic group for $i = s+1, \ldots, s+f$.
In any such decomposition the integer $f$ is uniquely determined and the prime-powers $p_i^{\alpha_i}$ are determined to within rearrangement.

Thus if $S$ is bad try factorizing $S$ to find the possible $p_i$. Then do Smith normal form calculations over $Z_{p_i^{\beta_i}}$ to find primary invariants rather than torsion invariants. Using $\beta_i$ we reveal all $p$-primary invariants with exponent less than $\beta_i$ explicitly. The choice of $\beta_i$ is made for convenience of computation. (After this we can assemble the primary invariants to find a more useful replacement for $S$.) If we want the isomorphism and $F$ is nontrivial we somehow have to factor out the torsionfree component and determine a homomorphism $G \to F$.

## 4 Lattice basis reduction

The rows (columns) of $A$ can be viewed as an integer lattice. Keeping entries small in computations with $A$ is closely related to finding small bases for integer lattices. Havas and Sterling [18] already used heuristic techniques for this, which proved to be time consuming but not fully effective. We have also implemented some other heuristic row reduction routines which have shown promising performance. For a row $\mathbf{r}$ of $A$, let $||\mathbf{r}||$ denote the sum of the absolute values of the entries in $\mathbf{r}$. Then one of our reduction routines is to consider every (ordered) pair $(\mathbf{r}, \mathbf{s})$ of distinct rows, and replace $\mathbf{r}$ by $\mathbf{r}' = \mathbf{r} + \mathbf{s}$ or $\mathbf{r}' = \mathbf{r} - \mathbf{s}$ if $||\mathbf{r}'|| < ||\mathbf{r}||$. We repeat this process until no further changes can be made. This routine is somewhat time consuming, but in many examples, if it is applied

just once or twice at critical times, then it can curtail and prevent further entry explosion. The difficulty is finding good heuristics to decide when to apply it. The corresponding routine can also be applied to columns of course, but it does not seem to be worthwhile doing both.

In their very important paper on computational number theory Lenstra, Lenstra, Lovász [25] included a new basis reduction algorithm which requires $O(X^4 \log(x))$ operations on numbers of length $O(X \log(x))$ and guarantees the quality of the reduced basis. The LLL algorithm was described for square matrices of full rank. Pohst [30] extended it to handle general rectangular matrices producing a modified algorithm, MLLL, with analogous complexity. MLLL produces a reduced basis from (possibly) linearly dependent vectors.

So, to get a homomorphism $G \to F$ we can use MLLL on the columns of $A$. (The infinite cyclic groups appear as columns of zeros.) In practice we start by doing some integer Gauss-Jordan eliminations first because eliminations are much faster than MLLL. We then factor $F$ out and compute both the structure of $T$ and a homomorphism $G \to T$ by modular methods. Here we have the usual type of compromise: the more integer Gauss-Jordan eliminations and less MLLL we do, the faster this works; however this also gives greater growth in intermediate and transforming matrix entries.

With a time efficiency perspective, since LLL is an $O(X^4)$ algorithm, it is worth reconsidering Gauss-Jordan elimination, which is $O(X^3)$ if we can keep down the size of entries.

## 5 Sparse matrices and pivoting strategies

In numerical analysis various techniques exist for handling sparse matrices. These have been much studied, and one recent reference is Zlatev [34]. Some applications of these kinds of techniques to exact matrices appear in Donald and Chang [8] and LaMacchia and Odlyzko [24], where some progress is made. In our context trying to find sparse initial relation matrices helps. Then use of pivoting strategies which both maintain sparsity and reduce entry explosion has given orders of magnitude improvement in performance. Observe that no particular pivoting strategy is specified in the algorithm given in §1: the algorithm simply says "choose a nonzero entry (the *pivot*)". We investigate various strategies in practice and find some which are not expensive to implement and which have good performance.

Note that Rose and Tarjan [31] have shown that it is an NP-complete problem to find a strategy which minimizes fill-in. This indicates that finding optimal pivoting strategies for our purposes is likely to be very difficult. Thus we concentrate on heuristics for obtaining good solutions.

In the $Z$-module context, it was rapidly observed that the pivot should be an element of minimal magnitude, and most implementations do this. For the kinds of modules we study, there are many unit ($\pm 1$) pivots. It turns out that careful selection among these is of great importance, and we use our pivoting strategies whenever there is more than one unit potential pivot.

To describe pivoting strategies we look at various row and column metrics related to potential pivots. For a vector $\mathbf{c} = (c_1, \ldots, c_n)$ we define: $|\mathbf{c}|_k = (\sum_1^n |c_i|^k)^{(1/k)}$ for $1 \leq k < \infty$; $|\mathbf{c}|_0 = \lim_{k \to 0}(\sum_1^n |c_i|^k)$; $|\mathbf{c}|_\infty = \lim_{k \to \infty} |\mathbf{c}|_k$. For $1 \leq k \leq \infty$ these metrics are the standard linear algebra norms, while for $k = 0$ the metric counts the number of nonzeros in the vector. Then, for $v = a_{s,t}$ a potential pivot in the $m \times n$ matrix $A$, we define

$$\begin{aligned}
|v|_k^R &= |(a_{s,1}, \ldots, a_{s,n})|_k \\
|v|_k^C &= |(a_{1,t}, \ldots, a_{m,t})|_k \\
|v|_k^{R-} &= |(a_{s,1}, \ldots, a_{s,t-1}, a_{s,t+1}, \ldots, a_{s,n})|_k \\
|v|_k^{C-} &= |(a_{1,t}, \ldots, a_{s-1,t}, a_{s+1,t}, \ldots, a_{m,t})|_k
\end{aligned}$$

Thus these are the corresponding metrics of the row and column containing $v$, possibly excluding $v$ itself.

For numerical applications Markowitz [26] introduced a heuristic where, in our terms, the pivot $v$ is chosen to minimize $|v|_0^{R-} \times |v|_0^{C-}$. Suitably modified for stability, this strategy remains a recommended one for general unsymmetric matrices (see Zlatev [26, §4.4]).

Our general idea is to consider various mathematical combinations of row and column metrics associated with potential pivots. We then choose pivots which minimize the combination of the associated row and column metric. Each combination can be interpreted as an estimate of some local property of the matrix. Thus, the Markowitz criterion clearly is an upper bound to the amount of fill-in which will occur if that pivot is used. It can be regarded as an approximation to the fill-in, and the ensuing strategy is thus an approximation to the greedy strategy with respect to minimizing fill-in.

Consider $|v|_1^{R-} \times |v|_1^{C-}$. For unit $v$, this product approximates the increase in the sum of the absolute values of all entries after pivoting on $v$. Or consider $|v|_\infty^{R-} \times |v|_\infty^{C-}$. For unit $v$, this product approximates the maximum increase in an entry after pivoting on $v$. It is very reasonable to look at minimizing these quantities for our purposes. Note that it is not particularly expensive to keep track of any one of these metrics for the whole matrix. After an initial computation, the metrics merely need updating for each row and column in which elements change.

There is another pivoting decision to be made, which is not specified in the simple algorithm description in §1. Thus, the algorithm reads: "While there is an entry $a_{1,j}$ in the first row not divisible by $a_{1,1}$: compute $a_{1,j} = a_{1,1}q + r$ ..." Once the chosen pivot is no longer a unit, it may well happen that there are many $a_{1,j}$ not divisible by the pivot (likewise for entries in the first column). In fact it usually is the case that there are many such nondivisible entries. (Given a random $a_{1,1}$, the chance that a random entry is divisible by it is $1/a_{1,1}$. As pivots increase, this probability decreases.)

So how should the $a_{1,j}$ (likewise $a_{j,1}$) be chosen? As before, many implementations simply take the first nondivisible entry found. As before, this is not a good idea.

What makes a nondivisible entry attractive for this purpose? If the remainder $r$ after division by $a_{1,1}$ is a unit ($\pm 1$ in the case of $Z$), then we can quickly proceed to the next step of the process. Basically, the smaller the magnitude and the number of divisors of the remainder, the better.

Thus, instead of choosing one entry $a_{1,j}$ and computing $a_{1,j} = a_{1,1}q + r$ for it, we do so for every other entry in the first row and column. (Currently we allow $r$ to be in the range $[-(|v|-1), |v|-1]$, though it is possible to halve the range.) This produces a complete row and column all of whose entries (excluding $a_{1,1}$ itself) are less than the previous minimum entry. This process can be seen as a special type of row and column reduction strategy, with the pivot entry as the key. In practice it has the same kind of beneficial effect on the matrix as other reductions and is very fast. We call this the best remainder strategy and use it in combination with our other pivoting strategies. The best remainder strategy combined with SGJ already makes significant improvements to performance. We observe that it often reintroduces a number of unit entries to the matrix, which are then subject to our other pivoting strategies.

# 6  Empirical results

## 6.1  Examples

Our examples come mainly from $Z$-modules which have arisen in actual group theoretic applications. Generally, the invariant factors in these cases are "nice". That is, there are a small number of finite nontrivial modules of small size. In these cases we start off with a matrix with small entries and finish with a diagonal matrix with small entries. Under these circumstances it seems reasonable to desire nice transforming matrices $P$ and $Q$. Since the initial and final matrices are identical regardless of algorithm, we measure the quality of solution in terms of the size of the biggest intermediate entry during the computation and of the size of the entries in the corresponding $P$ and $Q$. We have experimented with some artificial examples specially chosen in an effort to create poor performance for some algorithms. These produced behaviour consistent with the examples we present.

Our main examples come from module presentations which motivated the development of the Havas and Sterling implementation. The aim was to investigate certain groups by studying particular sections, namely abelian quotients of certain subgroups. This technique has led to effective understanding of some groups whose structure was not well-enough known, including the Fibonacci group $F(2,9)$ (see Havas, Richardson and Sterling [17] and Newman [29]) and certain knot groups (see Havas and Kovàcs [16]).

A common thread in these applications is the discovery of the module presentation by use of a Reidemeister-Schreier algorithm followed by abelianization. Analysis of this process reveals aspects of the nature of the ensuing presentation.

Consider the following situation. We start with a $d$-generator, $r$-relator presentation for a group $G$ and a subgroup $H$ of index $i$ whose abelian quotient $H/H'$ we wish to recognize. Then a direct implementation of the Reidemeister-Schreier method (see Havas [14]) leads to an integer relation matrix with $ri$ rows and $(d-1)i + 1$ columns. Most often $d$ and $r$ are small, but $i$ can be quite large, leading to sizable matrices.

The number of nonzero entries in the matrix can be readily bounded above. If the total relator length is $l$, then there can be no more than $(l/r) \times ri$ nonzeros, since each symbol in each relator converts to at most one nontrivial Schreier generator when rewritten. Even using this upper bound we see that the matrix will be sparse unless $l$ is relatively large, which is usually not the case. Thus the proportion of nonzero entries is bounded by $(li)/(ri \times ((d-1)i+1))$, which is approximately $l/((d-1)ri)$. This indicates why sparse matrix techniques are often applicable.

Of course each $m \times n$ matrix has $m!n!$ possibly different equivalent matrices under row and column permutations. We observe that different permutations do effect the performance of many of the integer algorithms, since some aspects of pivot selection rely on an arbitrary choice from among a number of equivalent potential pivots at each stage. This tends to be altered by order of entries in the matrix, with say the "first" of the equivalent pivots being selected.

## 6.2 Initial integer elimination methods

Consider the Fibonacci group $F(2,9) = \langle x_1, \ldots, x_9 \mid x_1 x_2 = x_3, \ldots, x_9 x_1 = x_2 \rangle$. The maximal nilpotent quotient of $F(2,9)$ is isomorphic to $Q_8 \times C_{19}$, a group of order 152. Hence we find subgroups of index 2, 4, 8, 19, 38, 76 and 152 in $F(2,9)$ corresponding with those of $Q_8 \times C_{19}$. The index 152 subgroup is interesting and we study its abelian quotient, which turns out to be $18C_5$.

To minimize the size of relation matrix it is natural to minimize the number of generators and relators for $F(2,9)$. We readily obtain a 2-generator, 2-relator presentation by eliminating redundant generators. Then the index 152 subgroup has a presentation with 304 relations and 153 generators, giving a $304 \times 153$ relation matrix for its abelian quotient. We denote by $A_1$ one specific such relation matrix.

Simple Gauss-Jordan (SGJ) elimination (choose as pivot the first element of minimal magni-tude) fails reasonably quickly. With 36 bit word size, overflow occurs after 104 eliminations. Using a multiple precision implementation, we find that entries with about 370 decimal digits exist after 120 eliminations.

Havas and Sterling [18] present a $26 \times 27$ matrix called $R_1$ which arises from a knot group. It makes an interesting test case: the initial matrix has rank 25, small entries, moderate density, and one nontrivial invariant, 3. (The matrix comes from a knot group presentation with 3 generators and 2 relators and a subgroup of index 13. The relators are unusually long, with total length 90. $R_1$ has 702 entries, 326 of which are zero, 281 of unit magnitude, 91 of magnitude 2, and 4 of magnitude 3.)

SGJ results in a maximum magnitude entry with 110 decimal digits. If we are less sensible and pivot on the first nonzero entry, then a maximum magnitude entry with 265 decimal digits occurs. Worse still, if we choose to pivot on a maximal magnitude entry each time, then a maximum magnitude entry with 1626 digits arises after just 12 eliminations. This strategy seems a potential candidate for proving exponential entry growth, but we have not succeeded in doing so.

We used versions of the algorithms of Kannan and Bachem [23] and of Chou and Collins [5] for Hermite normal form calculation. The Kannan-Bachem algorithm led to an 11 digit maximum magnitude entry (51665919764) in computing the HNF, from which the SNF is very easily calculated. Our experience shows that there is no appreciable difference between these algorithms for the matrices we have studied, but they are significantly better than SGJ.

When considered as a SNF problem for integer matrices there is no important difference due to matrix transposition. The result is simply transposed. When considered as an abelian group or module problem then there is a duality, with generators and relators interchanged. For the group situation it may be more convenient to emphasize operations on the relators rather than on the generators, because this makes an isomorphism computation simpler.

Notice that SNF algorithms generally imply significant differences in how often row operations are done, as against column operations. So it is worth reconsidering all results starting with the initial matrix transposed. Here we look at what happens starting with $R_1^T$.

SGJ results in a maximum magnitude entry with 1105 decimal digits, a factor of ten worse. Pivoting on the first nonzero entry, we obtain a maximum magnitude entry with 1050 decimal digits, which is (perhaps surprisingly) better than SGJ. Pivoting on a maximal entry each time, we obtain a maximum magnitude entry with 1704 digits after 13 eliminations, there having been 699 digits after 12 eliminations. In this case the maximum magnitude of the Kannan and Bachem algorithm was a 13 decimal digit number, 9330076432385.

## 6.3 Modular methods

The modular techniques described here are built into the computer algebra language *Cayley* [3]. When applied to finding the abelian quotient of the index 152 subgroup of $F(2,9)$ *Cayley* proceeds as follows.

*Cayley* uses primes with about 15 bits, so that all required computations can be conveniently done using 32 bit integer arithmetic. First a Hadamard bound is computed (in floating point) to determine how many primes are required to guarantee correct rank computation. Using the first prime (44657), *Cayley* discovers that the matrix has full (column) rank, so the Hadamard bound is not required for the rank computation. (While computing the rank the first modular determinant was also computed, at negligible extra cost.) The Hadamard bound indicates that at most 14 primes are required for determinant calculations (though there is also an early stopping criterion which may be used). Thus at worst 15 Gauss-Jordan eliminations are required so far.

For each prime used *Cayley* computes five modular determinants. The computation of five determinants rather than one takes little extra time, since the same calculations are done for the first (rank−1) rows, with five different rows, linearly independent from the first ones, being added as final row.

This in fact is a relatively difficult example. *Cayley* gives: determinant 1 = 36132812500000000; determinant 2 = 34004211425781250; determinant 3 = 26329040527343750; determinant 4 = 33142089843750000; determinant 5 = 55114746093750000.

The gcd of the first two determinants is 7629394531250, which is indeed the gcd of all five. Since this number does not fit in 32 bits *Cayley* proceeds to factorize it, getting $2 \times 5^{18}$. Then *Cayley* continues by attempting to compute the primary invariants, in this case requiring two more Gauss-Jordan eliminations.

For the prime 2 it computes modulo $2^2$, which is guaranteed to reveal any 2-factor exactly. It finds no 2-factor.

For the prime 5 it is less comfortable, since $5^{18}$ does not fit into 32 bits. It prefers to try single precision computation, so chooses to work modulo $5^{12} = 244140625$ which fits into a word. (Should there be a 5-factor with order exceeding $5^{11}$ this would reveal itself as zero in the SNF. This would necessitate going to multiple precision to determine its order.) In practice 18 five cycles appear.

As a result of this calculation it follows that $F(2,9)$ has a subgroup of index 190. We denote the abelian quotient of such a subgroup by $H_2$ and its relation matrix by $A_2$. $A_2$ has 380 rows and 191 columns, small entries and full rank.

The simple Gauss-Jordan elimination strategy gets to about $230 \times 40$ before overflow on a 32 bit machine. The modular method gives first determinant 1 770 749 945 013 406 400, second determinant 105 018 206 175 737 920. The gcd of these is $320 = 2^6 \times 5$, from which the associated abelian group is readily computed to be $C_4$.

After 90 SGJ eliminations on $A_2$ we obtain a $290 \times 101$ matrix with moderate sized entries. Applying the modular technique at this stage we obtain a first deteriminant with 55 digits:
6 809 695 809 169 251 595 747 179 546 442 846 834 847 729 901 468 148 960.
Such a number is hard to factorize. The second determinant also has 55 digits:
1 068 344 058 412 672 408 526 935 938 648 402 030 177 709 499 839 183 308. The gcd of these is 12, from which we can comfortably proceed.

A more spectacular example is provided by the Heineken group $G$ (see Neubüser and Sidki [28]), a group whose structure is still not well-enough understood. $G = \langle x, y, z \mid [x, [x, y]] = z, [y, [y, z]] = x, [z, [z, x]] = y \rangle$. We study sections to try to understand the group.

Consider $N = \langle a^5, b^5, (ab^{-2})^2, a^2baba^{-2}bab^{-1}a^{-1}b, aba^{-1}baba^{-1}bab^{-1}ab^{-1} \rangle_G$. Then $N$ has index 960 in $G$. We find the abelian quotient of $N$ using a Reidemeister-Schreier algorithm. Using a 2-generator presentation for G we get a relation matrix for $N/N'$ of size $1920 \times 961$ with full rank.

We use the notation $c_n$ for a composite number with $n$ decimal digits, $p_n$ for a prime with $n$ decimal digits, and $\hat{p}_n$ for a probable prime with $n$ digits. We obtain a Hadamard bound with some 3600 decimal digits. Then our first five determinants look like: $c_{175}$; $c_{174}$; $c_{175}$; $c_{175}$; $c_{175}$.

The gcd of determinants 1 and 2 is $2^{12} \times 3$, and modular diagonalization reveals that $N/N' \cong 8C_2 \oplus 2C_4$.

Note that factorization of numbers like these determinants (say $d_i$) is very hard. An elliptic curve method program of Richard Brent [1] reveals:
$d_1 = 2^{12}.3.7.11.71.109.1459.185914338563.c_{151}$;
$d_2 = 2^{12}.3.13.17.797.857.c_{162}$;
$d_3 = 2^{12}.3.11.89.52501.67153.303302071.283687.\hat{p}_{123}$;

$d_4 = 2^{12}.3.1233653.p_{14}.p_{21}.c_{131}$;
$d_5 = 2^{12}.3.181.284759.c_{163}$.

This shows the importance of taking the gcd of a few determinants when using the modular method. Remember that computing multiple determinants usually takes little extra time.

Even when the gcd has superficially seemed undesirable, factorization has shown otherwise. In one case a $576 \times 145$ matrix gave a 20 decimal digit gcd, which turned out to be $2^{65}$. The associated abelian group could then be readily computed as $4C_2$. A related $1152 \times 289$ matrix gave a 39 decimal digit gcd, which was $2^{122} \times 3^3$. This time the abelian group was $3C_2$.

## 6.4 Integer pivoting strategies

Consider $A_1$ ($304 \times 153$) of §6.2 again. As indicated, SGJ allowed 104 eliminations with a 36 bit word size before overflow. Using the simplest pivoting strategy (first nonzero entry) with multiple precision arithmetic led to entries with over 6000 decimal digits after 79 eliminations in one case. So 74 pivots were left.

| Pivoting strategies on permutations of $A_1$ ($304 \times 153$) | |
|---|---|
| Strategy | Performance |
| SGJ | all overflow with 37 to 50 pivots remaining |
| $R_0$ | all overflow with 22 to 41 pivots remaining |
| $C_0$ | all overflow with 33 to 41 pivots remaining |
| $\times_0$ | all succeed with maximum magnitude entry 44150 to 200415310 |
| $\times'_0$ | all succeed with maximum magnitude entry 44150 to 200415310 |
| $+_0$ | all succeed with maximum magnitude entry 3025 to 15586225 |
| $R_1$ | all overflow with 24 to 40 pivots remaining |
| $C_1$ | all overflow with 31 to 46 pivots remaining |
| $\times_1$ | all succeed with maximum magnitude entry 8120 to 264332655 |
| $\times'_1$ | all succeed with maximum magnitude entry 179645 to 264332655 |
| $+_1$ | all succeed with maximum magnitude entry 78230 to 21927855 |
| $R_2$ | all overflow with 15 to 42 pivots remaining |
| $C_2$ | all overflow with 33 to 44 pivots remaining |
| $\times_2$ | all succeed with maximum magnitude entry 182510 to 67404590 |
| $\times'_2$ | all succeed with maximum magnitude entry 1895 to 591771240 |
| $+_2$ | all succeed with maximum magnitude entry 3675 to 412676575 |
| $R_\infty$ | all overflow with 29 to 43 pivots remaining |
| $C_\infty$ | all overflow with 34 to 55 pivots remaining |
| $\times_\infty$ | 2 succeed with maximum magnitude entry 1925345 to 189210230 |
| | 6 overflow with 11 to 32 pivots remaining |
| $\times'_\infty$ | all overflow with 13 to 31 pivots remaining |
| $+_\infty$ | 2 succeed with maximum magnitude entry 214817492 to 617334795 |
| | 6 overflow with 20 to 26 pivots remaining |

The table shows the performance on 8 different permutations of $A_1$, in each case using the best remainder strategy. We designate the strategy which chooses a unit pivot $v$ with minimal $|v|_k^R$ by $R_k$, with minimal $|v|_k^C$ by $C_k$, with minimal $|v|_k^R + |v|_k^C$ by $+_k$, with minimal $|v|_k^R \times |v|_k^C$ by $\times_k$, and with minimal $|v|_k^{R-} \times |v|_k^{C-}$ by $\times'_k$.

Using SGJ combined with the best remainder strategy (now on a 32 bit machine), we found that all 8 overflow with from 37 to 50 pivots left. Looking at row metrics alone helps a little, but we do not complete in any case. All succeed with better pivoting strategies, taking both row and column metrics into account.

The table indicates clearly that combinations of the 0-, 1- and 2-metrics provide very effective pivoting strategies. Observe that the successful computations are much faster than the modular method, since only one matrix diagonalization is needed, compared with as many as 17 with the modular technique. Note that the metric based pivoting strategies without the best remainder strategy perform significantly worse. Maximum entries with 50 to 60 decimal digits arise. All this may be compared with the Kannan-Bachem method, which overflows 32 bit integers during the 77th elimination, and reaches a maximum magnitude entry with 35 decimal digits.

Given the variability shown here, we suggest that the best approach to such a problem is to start off by trying combinations of the 0-, 1- and 2-metrics. Only if unsuccessful that way should one revert to the modular method.

Similar results come from looking at $R_1$, where again we find these strategies dramatically outperforming the Kannan and Bachem method, which has proved polynomial complexity. The first four strategies (above the double line) do not incorporate the best remainder approach, while all the others do. Notice that the best remainder strategy alone added to SGJ makes it significantly outperform Kannan and Bachem.

| Pivoting strategies on $R_1$ ($26 \times 27$) ||
|---|---|
| Strategy | Maximum entry |
| Pivot on maximum | exceeds $10^{1626}$ after 12 eliminations |
| Pivot on first | $\sim 10^{265}$ |
| Pivot on minimum (SGJ) | $\sim 10^{110}$ |
| Kannan and Bachem | 51665919764 |
| SGJ + best $r$ | 152468 |
| $R_0$ | 17325 |
| $C_0$ | 141483 |
| $\times_0$ | 8501 |
| $\times_0'$ | 8501 |
| $+_0$ | 8501 |
| $R_1$ | 78928 |
| $C_1$ | 141483 |
| $\times_1$ | 57351 |
| $\times_1'$ | 8501 |
| $+_1$ | 32280 |
| $R_2$ | 46981 |
| $C_2$ | 22459 |
| $\times_2$ | 53613 |
| $\times_2'$ | 53613 |
| $+_2$ | 228030 |
| $R_\infty$ | 32210 |
| $C_\infty$ | 800722 |
| $\times_\infty$ | 41095 |
| $\times_\infty'$ | 34019 |
| $+_\infty$ | 66873 |

In this case it is easy to see how best remainder strategy has its impact. Best remainder does not come into play until the matrix has no unit entries. At that stage, after 11 pivoting steps on unit entries with SGJ, we have a diagonal of ones leading to a $14 \times 16$ working matrix, say A:

```
 -627  113   955  202 -282  164 -455 -139 -337  220 -114   85  261  473  113  240
  545  -95  -831 -176  246 -145  398  121  292 -193   99  -73 -230 -413  -96 -209
  179  -30  -269  -59   80  -46  130   40   94  -65   34  -23  -76 -135  -31  -69
  969 -167 -1471 -314  436 -258  706  215  517 -343  178 -129 -409 -733 -173 -371
  476  -87  -730 -152  214 -125  345  107  257 -167   85  -65 -197 -360  -86 -181
 -305   55   464   97 -138   80 -221  -67 -164  109  -55   41  127  231   53  116
 -271   49   417   87 -121   73 -197  -60 -147   96  -48   38  112  207   49  103
 -596  107   914  190 -269  157 -434 -133 -322  211 -105   81  247  452  105  228
  278  -52  -428  -88  122  -72  201   62  150  -99   48  -39 -112 -212  -50 -105
  435  -78  -664 -138  194 -114  315   97  232 -151   77  -59 -180 -327  -78 -167
 -218   39   332   71  -99   58 -158  -49 -117   78  -40   30   92  165   39   82
 -330   60   502  104 -148   86 -237  -73 -177  117  -59   45  136  249   59  125
 -521   94   799  168 -234  136 -380 -117 -283  183  -92   70  213  394   94  200
 -339   65   523  107 -153   87 -245  -77 -185  120  -59   47  137  257   61  128
```

The nonzero entry with minimal magnitude in the matrix is $a_{3,12} = -23$, and is used as a pivot. SGJ swaps it to the top left position, then computes the first nonzero remainder of an entry in the top row or leftmost column divided by the current pivot to become the next pivot. (There is a little flexibility about the details of these steps, which can lead to some variations in performance.) This process is repeated with the new pivot, and so on until we finally get a unit pivot, at which stage the row and column are zeroed. A simple implementation of this process gives a dense (no nonzero entries) $13 \times 15$ working matrix with minimal magnitude nonzero entry 3 and maximal magnitude entry 9968.

On the other hand, SGJ with a best remainder strategy computes a new $14 \times 16$ matrix first. (Our implementation does not move the minimal entry to the top left position at this stage, since it probably will not end up as a diagonal entry.) We obtain:

```
  22   7 -28  -7   6  -6  15  -3   9  -7   4  16 -15 -12   4 -15
 -20  -1  20   9  -6   1 -12  -3  -6  10  -7  -4  10  12   1  10
  18  -7 -16 -13  11   0  15  17   2 -19  11 -23  -7 -20  -8   0
 -24  -3  28   9  -6   0 -14   1  -9  10  -6 -14  13  12  -4  16
 -15  -8  17   4  -3   5 -10   8  -7   1  -2 -19  12   5  -5  14
   0   7  -3   2  -4  -2  -1  -9   2   8  -3  18  -3   6   4  -7
  13   4 -17  -2   4  -3   8  -5   7   1   1  15  -9  -3   3 -11
  25   5 -25 -11   7  -5  16  -1   8  -8   9  12 -17 -13   0 -15
 -13  -6  17   3  -6   6  -9   6  -8  -2  -2 -16  12   3  -3  12
 -14  -5  17   6  -5   4 -10   4  -8   5  -4 -13  11   8  -3  10
  10   2 -14  -2   2  -2   7  -2   5  -1   1   7  -5  -5   1  -8
   3   8  -9   1  -2  -4   3 -11   5   8  -3  22  -6   4   6 -10
  23   3 -19 -11   9  -4  15   4   3 -14  11   1 -18 -16   0 -10
  26   4 -26 -13  10  -7  20   4   7 -12  10   1 -18 -18  -2 -13
```

This is a nice example of very effective row and column reductions. Observe the appearance of 14 unit entries. Now it pivots on the first unit entry, in position (2,2), giving the following $13 \times 15$ matrix (which is clearly much better than that produced by SGJ):

$$\begin{array}{rrrrrrrrrrrrrr}
-118 & 112 & 56 & -36 & 1 & -69 & -24 & -33 & 63 & -45 & -12 & 55 & 72 & 11 & 55 \\
158 & -156 & -76 & 53 & -7 & 99 & 38 & 44 & -89 & 60 & 5 & -77 & -104 & -15 & -70 \\
36 & -32 & -18 & 12 & -3 & 22 & 10 & 9 & -20 & 15 & -2 & -17 & -24 & -7 & -14 \\
145 & -143 & -68 & 45 & -3 & 86 & 32 & 41 & -79 & 54 & 13 & -68 & -91 & -13 & -66 \\
-140 & 137 & 65 & -46 & 5 & -85 & -30 & -40 & 78 & -52 & -10 & 67 & 90 & 11 & 63 \\
-67 & 63 & 34 & -20 & 1 & -40 & -17 & -17 & 41 & -27 & -1 & 31 & 45 & 7 & 29 \\
-75 & 75 & 34 & -23 & 0 & -44 & -16 & -22 & 42 & -26 & -8 & 33 & 47 & 5 & 35 \\
107 & -103 & -51 & 30 & 0 & 63 & 24 & 28 & -62 & 40 & 8 & -48 & -69 & -9 & -48 \\
86 & -83 & -39 & 25 & -1 & 50 & 19 & 22 & -45 & 31 & 7 & -39 & -52 & -8 & -40 \\
-30 & 26 & 16 & -10 & 0 & -17 & -8 & -7 & 19 & -13 & -1 & 15 & 19 & 3 & 12 \\
-157 & 151 & 73 & -50 & 4 & -93 & -35 & -43 & 88 & -59 & -10 & 74 & 100 & 14 & 70 \\
-37 & 41 & 16 & -9 & -1 & -21 & -5 & -15 & 16 & -10 & -11 & 12 & 20 & 3 & 20 \\
-54 & 54 & 23 & -14 & -3 & -28 & -8 & -17 & 28 & -18 & -15 & 22 & 30 & 2 & 27 \\
\end{array}$$

Here we still have six potential pivots which are units left, in comparison to a minimal magnitude entry of 3 with straight SGJ. The maximal magnitude entry is 158, compared with 9968 for straight SGJ. Then this type of effect is magnified with further pivoting. The sequence of initial pivots (chosen when scanning for a nonzero entry in the working matrix) with the best remainder strategy is: 11 units, followed by -23, 1, -1, -12, 2, -1, -1, -24, 214, 760, 121, 7, 2754 and 6846. We have seen in detail what the row and column reductions do after the first nonunit pivot, -23. The next pivot -12 leads to a minimal entry of 2 in the ensuing matrix, and that leads to unit entries. Each of -24, 214 and 760 likewise lead to unit entries in two steps, the initial pivot 7 in one step and 2754 in four steps. Finally, initial pivot 6846 takes six steps to reach the entry 3, in a gcd computation for 6846 and -7551. With SGJ the sequence of initial pivots is: 11 units, followed by -23, 3, -26, -101747, -4366, -3320, -2330342737790992, -3619, -4171, -19612, -186162, -166597, -73125, and the last pivot is a 110 decimal digit number.

## 6.5 Sparsity considerations

Sparse matrix representations can reduce both the space required and the time required for matrix computations. Thus various algorithms which run in $O(X^3)$ time for dense matrices can be designed to run in $O(X^2)$ for sparse matrices which remain sparse. Hence the emphasis on finding pivoting strategies which minimize fill-in.

In §6.1 we wrote: "it is natural to minimize the number of generators and relators for $F(2,9)$". This turns out to be an oversimplification, based on the idea that minimizing the matrix dimensions is paramount. However, it is not just the matrix size which counts: sparsity matters, perhaps even more. This is especially true if a sensible sparse matrix representation is used.

There are obvious $n$-generator, $n$-relator presentations for $F(2,9)$ for $2 \leq n \leq 9$, produced by eliminating $(9-n)$ generators from the initial presentation. An $n$-generator, $n$-relator presentation leads to $(152(n-1)+1)$-generator, $(152n)$-relator presentations for the index 152 subgroup. Choosing $n=2$ minimizes the size of the relation matrix. However it also reduces the sparsity.

Our algorithms are included in the *quotpic* package (see Holt and Rees [19]). In that context we initially use a sparse matrix representation before converting to a standard array representation. During the sparse matrix phase we follow some Tietze transformation program principles (see Havas, Kenne, Richardson and Robertson [15]). Thus we perform short eliminations and (abelianized) substring searching till no further improvement is possible. The short eliminations correspond to well-selected pivoting operations (a pivot $v$ is used whenever a unit $v$ satisfying $|v|_1^{R-} \leq 1$ exists). The substring searching corresponds to row reduction heuristics, exemplified in the next subsection.

Consider the following results using SGJ. In each case 8 different permutations of the initial relation matrix were used (generated in natural group-theoretic ways). SGJ is applied after the Tietze transformations are done with the sparse representation. The initial relation matrices have dimensions $152n \times (152(n-1)+1)$, while the tabulated dimensions and densities are at the beginning of the SGJ phase. The initial density can be approximated by multiplying the given density by the ratio of the given matrix size to the initial size.

| $n$ | dimensions | nonzeros | density | performance |
|---|---|---|---|---|
| 2 | $300 \times 153$ | 2350 | 5.1% | all overflow with 40 to 52 pivots left |
| 3 | $446 \times 295$ | 2465 | 1.9% | all overflow with 37 to 51 pivots left |
| 4 | $553 \times 402$ | 2340 | 1.1% | all overflow with 30 to 47 pivots left |
| 5 | $601 \times 450$ | 2309 | 0.9% | all overflow with 32 to 45 pivots left |
| 6 | $653 \times 503$ | 2170 | 0.7% | all overflow with 19 to 30 pivots left |
| 7 | $693 \times 543$ | 2191 | 0.6% | 4 overflow with 19 to 25 pivots left, 4 succeed with maximum 11 177 145 to 81 407 370 |
| 8 | $711 \times 561$ | 2189 | 0.5% | all succeed with maximum 135 872 to 13 915 625 |
| 9 | $748 \times 598$ | 2244 | 0.5% | 2 overflow with 18 pivots left, 6 succeed with maximum 1 415 421 to 67 676 608 |

This clearly shows the benefit of increasing sparsity (decreasing density). It indicates how group theorists who wish to apply these methods should give careful consideration to the way in which they derive the presentations for their $Z$-modules.

## 6.6 Lattice basis reduction

Consider the Heineken group $G$ again, this time with $N$ the normal subgroup of index 120 with quotient $SL(2,5)$. Using the given 3-generator, 3-relator presentation for $G$ we obtain a $360 \times 241$ matrix for $N/N'$. The nontrivial invariants are $\{2,2,2,2\}$. Using the $\times_1'$ pivoting strategy together with the best remainder strategy gives integer overflow with 8 pivots left. Use of the same methods, together with the row reduction routine parameterized to come into play the first time a matrix entry exceeds 1000 (and then 2000, 4000, etc.), results in just one application of row reduction when the working matrix is $135 \times 16$. This leads to successful completion with maximum magnitude entry 1319. The total execution time is less than double that required to reach the previous integer overflow.

Some details of this calculation are as follows. The $127 \times 8$ matrix obtained prior to overflow using the $\times_1'$ pivoting strategy together with the best remainder strategy is completely dense (no nonzero entries). Its first five rows were:

```
-3838608  2675947 -2212100   323972  2163968  4944023 -16113986 -506210
 4157363 -2897967  2394459  -350613 -2343527 -5354460  17451974  549329
 4443281 -3097387  2560331  -375045 -2504759 -5722740  18651742  586177
-2285021  1592690 -1315796   192806  1287914  2943010  -9591208 -302339
 -903675   630414  -523497    76793   509785  1163755  -3793754 -116721
```

Eight pivoting steps earlier, when the matrix was $135 \times 16$, the first five rows were:

```
   -49  147   46 -138  164  116 -224  -48  421 -160  330 -128  120 -471  -76 -242
     3    2   -1    5  -12   -9    1   -2  -31    8  -45   -3   15   48  -17   38
     2   -4   -1   -1   -1   -1    7    1   -3   -1   -8    5   -5   10    8    1
    29 -198    3  184 -124  -57  188   39 -245  150   91  142 -238   64  149   24
   -78  104  -50  -15   88  -58  116  127  467 -180  611   49 -237 -632  592 -631
```

The row including the largest entry, 1309, was:

197 -433 73 70 -201 173 -266 -283 -1014 401 -1090 -26 466 1233 -1198 1309.

The matrix was pretty dense (36 zero entries out of 2160), with 43 unit entries still available as potential pivots. Row simplification reduced the matrix to $73 \times 16$, since 62 linearly dependent rows were found in the process. The first five rows were:

```
0  0 0 -1 0  0 0 0 0 -1 0 0 0 0  0 0
0 -1 0  1 0  0 0 0 0  0 0 0 0 0  0 0
0  0 0 -1 0  1 0 0 0  0 0 0 0 0  0 0
1  0 1  0 1  0 1 0 0  0 0 0 0 0  0 0
0  0 0  0 1 -1 0 0 0  0 0 1 0 0 -1 0
```

The largest magnitude entries were $\pm 2$, all 15 of which occurred in rows with no other nonzero entry. The ensuing matrix with eight columns (where previously the overflow occurred) was $63 \times 8$ with similar structure. The largest magnitude entries were $\pm 2$ as before, and the first five rows were:

```
 0 0  0  0 0  1  0 -1
 1 1 -1  1 0  0  0  0
 1 0  0  1 0 -1 -1  0
 1 1  1 -1 0  0  0  0
-1 1 -1  1 0  0  0  0
```

The row reductions we apply this way are done according to heuristics. An alternative approach, which guarantees the quality of the reduced basis in a certain way, is to use the MLLL algorithm. However, even with MLLL, we have to decide when to apply it. Since MLLL has $O(X^4)$ complexity it is relatively expensive, and very slow with large matrices.

Thus, in both cases, we need some rule for initiating the row reduction process. Natural rules include: as soon as the matrix has no more unit entries; as soon as matrix entries exceed a certain size. Again these are heuristics, and our implementations include a number of them.

Consider a situation in which we want to know the canonical form of the group and an isomorphism. This means that we want the unimodular matrix $Q$, and we usually want it to have small entries.

Our solutions involve using sensible pivoting plus row reduction. Consider the $304 \times 153$ matrix $A_1$ again. With moderate use of MLLL we can obtain the SNF with maximum magnitude entry of 385 in the working matrix, maximum magnitude entry 145 in $Q$, and (for the record) maximum magnitude entry 9243 in $P$. Heavier (and much more time consuming) use of MLLL leads to a maximum magnitude entry of 101 in the working matrix, maximum magnitude entry 90 in $Q$, and maximum magnitude entry 572 in $P$.

# 7 Conclusions

We have described previous and new methods for recognizing badly presented $Z$-modules. The new methods are controlled by heuristics, but show dramatically improved performance over previous methods which have polynomial complexity. Empirical evidence of this is based mainly on presentations arising from group theoretic calculations. Since finding the best solution to these problems is a very difficult we suggest that heuristic methods like these are likely to provide the best practical solutions. The nature of the heuristics leaves us without formal complexity results, merely the empirical evidence.

The heuristics that we describe perform particularly well. As with all hard problems which are solved by heuristics, we recommend trying different heuristics when attempting to solve a problem which is not immediately resolved. The key heuristics are good pivot selection, together with lattice basis reduction if better solutions are required or more difficult problems are being solved.

# 8 Acknowledgements

We are grateful to Dr Keith Matthews of the University of Queensland who provided computational assistance, especially with infinite precision computations and the MLLL algorithm.

# References

[1] R.P. Brent, Some integer factorization algorithms using elliptic curves, *Australian Computer Science Communications* 8: 149–163 (1986).

[2] S. Cabay and T.P. Lam, Congruence techniques for the exact solution of integer systems of linear equations, *ACM Trans. Math. Software* 3: 386–397 (1977).

[3] J. Cannon and W. Bosma, *A handbook of Cayley functions*, Computer Algebra Group, University of Sydney, 1991.

[4] J. Cannon and G. Havas, Algorithms for groups, *Austral. Comp. J.* 24: 51–60 (1992).

[5] T.-W.J. Chou and G.E. Collins, Algorithms for the solution of systems of linear Diophantine equations, *SIAM J. Comput.* 11: 687–708 (1982).

[6] P.D. Domich, Residual Hermite normal form computations, *ACM Trans. Math. Software* 15: 275–286 (1989).

[7] P.D. Domich, R. Kannan and L.E. Trotter Jr., Hermite normal form computation using modulo determinant arithmetic, *Math. Oper. Res.* 12: 50–59 (1987).

[8] B.R. Donald and D.R. Chang, On the complexity of computing the homology type of a triangulation, in *Proc. 32nd Annual Symposium on Foundations of Computer Science*, 650–661, IEEE Computer Society Press, 1991.

[9] M.A. Frumkin, An application of modular arithmetic to the construction of algorithms for solving systems of linear equations, *Soviet Math. Dokl.* 17: 1165–1168 (1976).

[10] M.A. Frumkin, Polynomial time algorithms in the theory of Linear Diophantine equations, in *Fundamentals of computation theory*, Lecture Notes in Computer Science 56, Springer-Verlag, 1977, pp. 386–392.

[11] L.J. Gerstein, A local approach to matrix equivalence, *Linear Algebra Appl.* 16: 221–232 (1977).

[12] J.L. Hafner and K.S. McCurley, Asymptotically fast triangularization of matrices over rings, *SIAM J. Comput.* 20: 1068–1083 (1991).

[13] B. Hartley and T.O. Hawkes, *Rings, modules and linear algebra*, Chapman and Hall, London, 1970.

[14] G. Havas, A Reidemeister-Schreier program, in *Proc. Second Internat. Conf. Theory of Groups*, Lecture Notes in Mathematics 372, Springer-Verlag, Berlin, Heidelberg, New York, 1974, pp. 347–356.

[15] G. Havas, P.E. Kenne, J.S. Richardson and E.F. Robertson, A Tietze transformation program, in *Computational Group Theory* (M.D. Atkinson, ed.), Academic Press, London, 1984, pp. 67–71.

[16] G. Havas and L.G. Kovàcs, Distinguishing eleven crossing knots, in *Computational Group Theory* (M.D. Atkinson, ed.), Academic Press, London, 1984, pp. 367–373.

[17] G. Havas, J.S. Richardson and L.S. Sterling, The last of the Fibonacci groups, *Proc. Roy. Soc. Edinburgh*, 83A: 199–203 (1979).

[18] G. Havas and L.S. Sterling, Integer matrices and abelian groups, in *Symbolic and Algebraic Computation*, Lecture Notes in Computer Science 72, Springer-Verlag, 1979, pp. 431–451.

[19] D.F. Holt and S. Rees, A graphics system for displaying finite quotients of finitely presented groups, in *Proc. DIMACS Workshop Groups and Computation, 1991*, AMS-ACM, to appear.

[20] T.C. Hu, *Integer programming and network flows*, Addison Wesley, 1969.

[21] C.S. Iliopoulos, Worst case complexity bounds on algorithms for computing the canonical structure of finite abelian groups and the Hermite and Smith normal forms of an integer matrix, *SIAM J. Comput.* 18: 658–669 (1989).

[22] C.S. Iliopoulos, Worst case complexity bounds on algorithms for computing the canonical structure of infinite abelian groups and solving systems linear Diophantine equations, *SIAM J. Comput.* 18: 670–678 (1989).

[23] R. Kannan and R. Bachem, Polynomial Time algorithms for computing Smith and Hermite normal forms of an integer matrix, *SIAM J. Comput.* 8: 499–507 (1979).

[24] B.A. LaMacchia and A.M. Odlyzko, Solving large sparse linear systems over finite fields, *Advances in cryptology — CRYPTO '90*, Lecture Notes in Computer Science 537, Springer-Verlag, 1990, pp. 109–133.

[25] A.K. Lenstra, H.W. Lenstra Jr. and L. Lovász, Factoring polynomials with rational coefficients, *Math. Ann.* 261: 515–534 (1982).

[26] H.M. Markowitz, The elimination form of the inverse and its application to linear programming, *Management Sci* 3: 255–269 (1957).

[27] J. Neubüser, An elementary introduction to coset-table methods in computational group theory, *Groups — St Andrews 1981*, London Math. Soc. Lecture Note Ser. 71, Cambridge University Press, Cambridge, 1984, pp. 1–45.

[28] J. Neubüser and S. Sidki, Some computational approaches to groups given by finite presentations, *Mathematica Universitaria* 7: 77–120 (Rio de Janeiro, 1988).

[29] M.F. Newman, Proving a group infinite, *Arch. Math.* 54: 209–211 (1990).

[30] M. Pohst, A modification of the LLL reduction algorithm, *J. Symbolic Computation* 4: 123–127 (1987).

[31] D.J. Rose and R.E. Tarjan, Algorithmic aspects of vertex elimination directed graphs, *SIAM J. Appl. Math.* 34: 176-197 (1978).

[32] C.C. Sims, *Computation with finitely presented groups*, Cambridge University Press, Cambridge, to appear.

[33] H.J.S. Smith, On systems of linear indeterminate equations and congruences, *Philos. Trans. Royal Soc. London* cli: 293–326 (1861). See also: *The Collected Mathematical Papers of Henry John Stephen Smith*, Volume 1, New York: Chelsea (1965).

[34] Z. Zlatev, *Computational methods for general sparse matrices*, Kluwer Academic Publishers, Dordrecht, 1991.